\documentclass{article}
    \usepackage[parfill]{parskip}
    \usepackage[utf8]{inputenc}
    \usepackage{amsmath,amssymb,amsfonts,amsthm}
    
    \usepackage[english]{babel}

    \newtheorem{theorem}{Theorem}
    \newtheorem{corollary}{Corollary}

    \usepackage{graphicx,csquotes}

    \usepackage{authblk}

    \title{A $q$-analog of Jacobi's two squares formula and its applications}
    \author[1]{Jos\'e Manuel Rodr\'iguez Caballero}
    \affil[1]{D\'epartment de math\'ematiques et de statistique, Universit\'e Laval, Qu\'ebec, Canada.}
    \affil[ ]{\textit{jose-manuel.rodriguez-caballero.1@ulaval.ca}}
    \date{\today}

\providecommand{\keywords}[1]
{
  \small	
  \textbf{\textbf{Keywords:}} #1
}

\providecommand{\msc}[1]
{
  \small	
  \textbf{\textbf{Mathematics Subject Classification:}} #1
}
    
\begin{document}

\maketitle

\begin{abstract}
    We consider a $q$-analog $r_2(n, q)$ of the number of
representations of an integer as a sum of two squares $r_2(n)$. This
$q$-analog is generated by
the expansion of a product that was studied by Kronecker and Jordan. We generalize
Jacobi's two squares formula from $r_2(n)$ to $r_2(n, q)$. We characterize the
signs in the coefficients of $r_2(n, q)$ using the prime factors of $n$. 
We use $r_2(n, q)$ to characterize the integers which are the length of the hypotenuse
of a primitive Pythagorean triangle.
\end{abstract}


\msc{11C08, 11E25, 11N13.}

\section{Introduction}

Using the formal identity 
\begin{equation}\label{djkljelkfjwewj4384rsdsfsf}
    \prod_{m=1}^{\infty} \frac{1 -(-t)^m}{1 + (-t)^m} = \sum_{n=-\infty}^{\infty}  t^{n^2},     
\end{equation}
due to C. F. Gauss \cite{gauss1866werke}, the generating
function for the number of representations of an integer $n$ as the sum of the squares of two integers, denoted $r_2(n)$,
immediately follows,

\begin{equation}\label{rwrsfsfht7i8DFGHty76UT}
    \prod_{m=1}^{\infty} \frac{\left( 1 -(-t)^m \right)^2}{\left( 1 + (-t)^m \right)^2} = 1+\sum_{n=1}^{\infty}  r_2(n) t^n.
\end{equation}

C. G. J. Jacobi \cite{jacobi1829fundamenta} expressed $r_2(n)$ as a function of the divisors of $n \geq 1$,
\begin{equation}\label{tetERetgdfgd345343}
    r_2\left(n\right) = 4 d_{1,4}(n) - 4 d_{3,4}(n), 
\end{equation}
where $d_{k,m}(n)$ is the number of divisors of $n$ which are congruent to $k$ modulo $m$.

The explicit formula for $r_2(n,q)$ defined by the expansion of the product

\begin{equation}\label{rrwRWEewrwsdFSDFxc5464}
    \prod_{m=1}^{\infty} \frac{\left( 1 -(-qt)^m \right)^2}{\left(1 + q(-qt)^m \right)\left(1 + q^{-1}(-qt)^m \right)}
= 1+\sum_{n=1}^{\infty} r_2\left(n, q \right) t^n,
\end{equation}
which is a $q$-deformation\footnote{A $q$-deformation of an expression $A(t)$ is another expression $A(t,q)$
satisfying $A(t,1) = A(t)$.} of identity \eqref{rwrsfsfht7i8DFGHty76UT},
can be attributed to L. Kronecker \cite{kronecker1890theorie}, who proved a more general version
 and C. Jordan \cite{jordan1894fonction}, who described a method to derive this particular case. 
 The polynomial $r_2\left(n, q \right)$
was first introduced in the article \cite{rodriguez2019characterization}, using the notation $\Gamma_n(q)$, and it was
called \emph{Kassel–Reutenauer $q$-analog of the number of representations as a sum of two squares}. Some versions of the 
polynomial $r_2\left(n, q \right)$, e.g., changing the sign of $q$ and sometimes dividing it by $q-1$ or by $(q-1)^2$, have been studied by several 
 authors, in connection with different branches of mathematics: finite fields \cite{kassel2018complete, kassel2018counting, caballero2018kassel},
algebraic topology \cite{hausel2013arithmetic}, modular functions \cite{kassel2017fourier}
and elementary number theory \cite{caballero2019function, caballero2020jordan, caballero2019integers}.

The aim of the present note is to prove that there are polynomials $4d_{1,4}(n,q)$ and $4d_{3,4}(n,q)$ satisfying the following properties\footnote{We consider that
it is more elegant to work with $4d_{1,4}(n,q)$ and $4d_{3,4}(n,q)$ rather than $d_{1,4}(n,q)$ and $d_{3,4}(n,q)$ because of property (i).}:

(i) all the coefficients of $4d_{1,4}(n,q)$ and $4d_{3,4}(n,q)$ are non-negative integers;

(ii) for every $k$,  $q^k$ cannot appear with non-zero coefficient in both $4d_{1,4}(n,q)$ and $4d_{3,4}(n,q)$;

(iii) the decomposition of $r_2(n, q)$ into positive and negative parts,
\begin{equation}\label{a3546456TRETERdfgdgdgd}
    r_2(n, q) = 4d_{1,4}(n,q) - 4d_{3,4}(n,q)
\end{equation}
holds;

(iv) $4d_{1,4}(n,q)$ and $4d_{3,4}(n,q)$ are $q$-analogs of $4d_{1,4}(n)$ and $4d_{3,4}(n)$ respectively, i.e.,
$4d_{1,4}(n,1) = 4d_{1,4}(n)$ and $4d_{3,4}(n,1) = 4d_{3,4}(n)$.

Therefore, the identity \eqref{a3546456TRETERdfgdgdgd}
is a generalization of Jacobi's two squares formula from integers to polynomials. Considering that properties (i), (ii)
and (iii) uniquely define the polynomials $4d_{1,4}(n,q)$ and $4d_{3,4}(n,q)$, it is non-trivial that they should also satisfy
property (iv). Furthermore,
as applications of this formula, we will determine when $r_2(n, q)$ has a negative coefficient by analyzing the prime factors of $n$.
Also, we will use $r_2(n, q)$ to characterize the integers which are the length of the hypotenuse of a primitive Pythagorean triangle.

\section{Generalization of Jacobi's formula}
In this section we will prove our main result.

\begin{theorem}\label{MainTheo}
    Let $n$ be a positive integer. The polynomials
    \begin{eqnarray}
        \label{sfkjsflkj34u43u0wekofp} 4d_{1,4}(n, q)
        &=&(q+1)\sum_{\scalebox{0.7}{$\begin{array}{c}
            d | n \\
            d \equiv 1 \pmod{4} 
        \end{array}$}} \left( q^{\left(\frac{2n}{d}+1\right)(d-1)/2} + q^{\left(\frac{2n}{d}-1\right)(d+1)/2}  \right), \\
       \label{fsdfvxvt632423dHGFH0} 4d_{3,4}(n, q)
        &=&(q+1)\sum_{\scalebox{0.7}{$\begin{array}{c}
            d | n \\
            d \equiv 3 \pmod{4} 
        \end{array}$}} \left( q^{\left(\frac{2n}{d}+1\right)(d-1)/2} + q^{\left(\frac{2n}{d}-1\right)(d+1)/2}  \right),
    \end{eqnarray} 
    satisfy properties (i), (ii), (iii) and (iv).
\end{theorem}

\begin{proof}

    Property (i) immediately follows from the explicit expressions \eqref{sfkjsflkj34u43u0wekofp} and \eqref{fsdfvxvt632423dHGFH0}. Property
    (iv) is just the result of the evaluations of these expressions at $q = 1$.

    Property (iii) follows from the following formal manipulation. Take formula (0.100) from \cite{cooper2017ramanujan},

    \begin{equation}\label{XVvxvvnnvhjkyu67543}
        \prod_{m=1}^{\infty} \frac{\left( 1 -t^m \right)^2}{\left(1 - qt^m \right)\left(1 - q^{-1}t^m \right)}
    = 1 + \left(q^{1/2} - q^{-1/2}\right) \sum_{d=1}^{\infty}\sum_{k=1}^{\infty} t^{d k} \left( q^{d-k/2} - q^{-d+k/2} \right).
    \end{equation}

    Replace $t$ by $qt$ in \eqref{XVvxvvnnvhjkyu67543},

    \begin{equation}\label{xcvxvTERT4535GDGDjgiiyiHK34}
        \prod_{m=1}^{\infty} \frac{\left( 1 -(qt)^m \right)^2}{\left(1 - q(qt)^m \right)\left(1 - q^{-1}(qt)^m \right)} 
        = 1 + \left(q - 1\right) \sum_{d=1}^{\infty}\sum_{ \scalebox{0.7}{$\begin{array}{c}
        k \geq 1 \\
        k \text{ odd}
        \end{array}$} } t^{d k} \left( q^{dk+ d -k/2 - 1/2} - q^{dk -d + k/2 - 1/2} \right),
    \end{equation}
where the restriction to only odd values of $k$ is because of the identity

\begin{eqnarray}
    & & \sum_{d=1}^{\infty}\sum_{ \scalebox{0.7}{$\begin{array}{c}
        k \geq 1 \\
        k \text{ even}
        \end{array}$} } t^{d k} \left( q^{dk+ d -k/2 - 1/2} - q^{dk -d + k/2 - 1/2} \right) \\
        &=& \sum_{d=1}^{\infty}\sum_{e = 1}^{\infty} t^{2 d e} \left( q^{2de+ d - e - 1/2} - q^{2de - d + e - 1/2} \right) \\
        &=& q^{-1/2} \left( \sum_{d=1}^{\infty}\sum_{e = 1}^{\infty} t^{2 d e}q^{2de + d - e} - \sum_{d=1}^{\infty}\sum_{e = 1}^{\infty} t^{2 d e}  q^{2de - d + e} \right) \\
        &=& 0.
\end{eqnarray}

Replace $q$ by $-q$ in \eqref{xcvxvTERT4535GDGDjgiiyiHK34},
    \begin{eqnarray}\label{xvxvadaADerwrERW45535}
        & & \prod_{m=1}^{\infty} \frac{\left( 1 -(-qt)^m \right)^2}{\left(1 + q(-qt)^m \right)\left(1 + q^{-1}(-qt)^m \right)} \\
        &=& 1 + \left(q + 1\right) \sum_{d=1}^{\infty}\sum_{ \scalebox{0.7}{$\begin{array}{c}
            k \geq 1 \\
            k \text{ odd}
            \end{array}$} } t^{d k} \left( (-q)^{(2d+1)(k-1)/2} - (-q)^{(2d-1)(k+1)/2}  \right) \\
        &=&  1 + \left(q + 1\right) \sum_{d=1}^{\infty}\sum_{ \scalebox{0.7}{$\begin{array}{c}
            k \geq 1 \\
            k \equiv 1 \pmod{4}
            \end{array}$} } t^{d k} \left( (-q)^{(2d+1)(k-1)/2} - (-q)^{(2d-1)(k+1)/2}  \right) \nonumber \\
         & &  + \left(q + 1\right) \sum_{d=1}^{\infty}\sum_{ \scalebox{0.7}{$\begin{array}{c}
                k \geq 1 \\
                k \equiv 3 \pmod{4}
            \end{array}$} } t^{d k} \left( (-q)^{(2d+1)(k-1)/2} - (-q)^{(2d-1)(k+1)/2}  \right) \\
            &=&  1 + \left(q + 1\right) \sum_{d=1}^{\infty}\sum_{ \scalebox{0.7}{$\begin{array}{c}
                k \geq 1 \\
                k \equiv 1 \pmod{4}
                \end{array}$} } t^{d k} \left( q^{(2d+1)(k-1)/2} + q^{(2d-1)(k+1)/2}  \right) \nonumber \\
             & &  - \left(q + 1\right) \sum_{d=1}^{\infty}\sum_{ \scalebox{0.7}{$\begin{array}{c}
                    k \geq 1 \\
                    k \equiv 3 \pmod{4}
                \end{array}$} } t^{d k} \left( q^{(2d+1)(k-1)/2} + q^{(2d-1)(k+1)/2}  \right) \\
            &=& 1 + \sum_{n=1}^{\infty} \left( 4d_{1,4}(n,q) - 4d_{3,4}(n,q) \right) t^n.
    \end{eqnarray}

    To prove property (ii), we proceed by \emph{reductio ad absurdum}. Suppose that for some $k$, the coefficient of $q^k$ is 
    non-zero in both $4d_{1,4}(n,q)$ and $4d_{3,4}(n,q)$. We need to analyze 16 possible cases. We will use the notations $d$ 
    and $e$ for two arbitrary divisors of $n$
    satisfying $d\equiv 1\pmod{4}$ and $e\equiv 3\pmod{4}$.

    Let $f(x) =  \left(\frac{2n}{x}+1\right)\frac{x-1}{2}$. Notice that $f(x)$, on the domain $x > 0$,
    is strictly increasing and satisfies the inequality $f(x+2) - f(x) > 1$. Furthermore,
    $f\left( \frac{2n}{x} \right) = \left(\frac{2n}{x}-1\right)\frac{x+1}{2}$.
    
    Notice that, $k = f(d) = f(e)$ implies $d = e$. Nevertheless, this is impossible because $d \not\equiv e \pmod{4}$. In the same way,
    it is easy to prove that $k = f\left(\frac{2n}{d}\right) = f\left(\frac{2n}{d}\right)$ also implies an absurde.
    Similarly, $k = f(d)+1 = f(e)+1$ and $k = f\left(\frac{2n}{d}\right)+1 = f\left(\frac{2n}{d}\right)+1$ are impossible.

    Notice that $k = f\left(d\right) = f\left(\frac{2n}{e}\right)$ implies $d = \frac{2n}{e}$. If this is the case, $d e = 2n$, which is 
    absurd, since $d$ and $e$ are odd. In the same vein, $k = f\left(\frac{2n}{d}\right) = f\left(e\right)$ is also absurd. Similarly,
    we exclude the cases $k = f\left(d\right)+1 = f\left(\frac{2n}{e}\right)+1$ and $k = f\left(\frac{2n}{d}\right)+1 = f\left(e\right)+1$.

    Assume that $k = f\left(d\right)+1 = f\left(e\right)$. It follows that $e > d$. Because $e$ and $d$ share the same parity (both are odd),
    $e \geq d + 2$. Hence, $f(e) - f(d) \geq f(d+2) - f(d) > 1$, which contradicts our assumption. In the same vein, it is easy prove that
    $k = f\left(d\right) = f\left(e\right)+1$ implies an absurd conclusion. Similarly, we exclude the cases,
    $k = f\left( \frac{2n}{d}  \right)+1 = f\left( \frac{2n}{e}\right)$ and $k = f\left( \frac{2n}{d}  \right) = f\left( \frac{2n}{e}\right)+1$
    by considering that $\frac{2n}{d}$ and $\frac{2n}{e}$ share the same parity (both are even).

    Assume that $k = f\left(d\right)+1 = f\left(\frac{2n}{e}\right)$. On the one hand, $f(d) =  \left(\frac{2n}{d}+1\right)\frac{d-1}{2}$ is even,
    since $\frac{d-1}{2}$ is even. On the other hand $f\left(\frac{2n}{e}\right) = \left(\frac{2n}{e}-1\right)\frac{e+1}{2}$ is also even, since $\frac{e+1}{2}$
    is even. We derive the absurd conclusion that $1 = f\left(\frac{2n}{e}\right) - f\left(d\right)$ should be even. In the same vein,
    we can easily prove that $k = f\left(d\right) = f\left(\frac{2n}{e}\right)+1$ imples an absurd. Similarly, we can exclude the cases
    $k = f\left(\frac{2n}{d}\right)+1 = f\left(e\right)$ and $k = f\left(\frac{2n}{d}\right) = f\left(e\right)+1$.

\end{proof}

\section{Applications}

In this section we derive some immediate consequences of our generalization of Jacobi's formula.

\begin{corollary} \label{coro1}
    Let $n$ be a positive integer. The polynomial $r_2\left(n, q \right)$ has a negative coefficient if and only if some of the prime factors of
    $n$ are congruent to $3$ modulo $4$.
\end{corollary}

\begin{proof}
    Considering that $d_{3,4}(n,q) \neq 0$ if and only if some of the prime factors of
$n$ are congruent to $3$ modulo $4$, the result immediately follows from Theorem \ref{MainTheo} and the definition of $d_{1,4}(n,q)$ and $d_{3,4}(n,q)$.

\end{proof}

We recall that $n$ is the hypotenuse of a primitive Pythagorean triangle 
if and only if for some pair of positive integers $u$ and $v$ the equality $u^2 + v^2 = n^2$ holds and $u$, $v$ and $n$ are relatively prime.

\begin{corollary}\label{coro2}
    An odd integer $n$ larger than $1$ is the length of the hypotenuse of a primitive Pythagorean
    triangle if and only if all the coefficients of 
    the polynomial $r_2\left(n, q \right)$ are non-negative.
\end{corollary}

\begin{proof}
    E. J. Eckert \cite{eckert1984group} proved that an integer larger than $1$ is the hypotenuse of a primitive Pythagorean triangle if and only if all its 
prime factors are congruent to $1$ modulo $4$. Combining this result with Corollary \ref{coro1}, the result follows.
\end{proof}

\section{Final remarks}

In the spirit of the work of C. Kassel and C. Reutenauer \cite{kassel2018counting}, the value of the polynomial $r_2(n, q)$,
when $q$ is a prime power, may have a combinatorial interpretation in the ring $\mathbb{F}_q[X,Y,X^{-1},Y^{-1}]$.

Let $r_4(n)$ be the number of representations of $n$ as the sum of $4$ squares of integers.
We suggest to empirically study the $q$-analog of $r_4(n)$ obtained from 
the square
\begin{equation}\label{fsfsdjlkjl534535werRWER}
    \left(1 + \sum_{n=1}^{\infty} r_2\left(n, q \right) t^n \right)^2 = 1 + \sum_{n=1}^{\infty} r_4\left(n, q \right) t^n
\end{equation}
and check whether some of the classical results about $r_4(n)$ can be generalized to $r_4(n, q)$.
The expansion of the corresponding product can be found in equation (0.101) of \cite{cooper2017ramanujan}.

\end{document}